\documentclass{article}
\usepackage{quad}

\hypersetup{pdftitle={5-point CAT(0) spaces after  Tetsu Toyoda},
pdfauthor={Nina Lebedeva and Anton Petrunin}}

\newcommand{\Addresses}{{\bigskip\footnotesize

\noindent Nina Lebedeva,
\par\nopagebreak
 \textsc{Saint Petersburg State University, 7/9 Universitetskaya nab., St. Petersburg, 199034, Russia}
\par
\nopagebreak
 \textsc{St. Petersburg Department of V.A. Steklov Institute of Mathematics of the Russian Academy of Sciences, 27 Fontanka nab., St. Petersburg, 191023, Russia}
  \par\nopagebreak
  \textit{Email}: \texttt{lebed@pdmi.ras.ru}

\medskip

\noindent   Anton Petrunin, 
\par\nopagebreak
 \textsc{Math. Dept. PSU, University Park, PA 16802, USA.}
  \par\nopagebreak
  \textit{Email}: \texttt{petrunin@math.psu.edu}
  
}}

\begin{document}

\title{5-point CAT(0) spaces after Tetsu Toyoda}
\author{Nina Lebedeva and Anton Petrunin}
\date{}
\maketitle
\begin{abstract}
We give another proof of Toyoda's theorem that describes 5-point subspaces in $\CAT(0)$ length spaces.
\end{abstract}

\nofootnote{\textbf{Keywords:} CAT(0), finite metric space, comparison inequality, Alexandrov comparison.}
\nofootnote{\textbf{MSC:} 53C23, 30L15, 51F99.}

\section{Introduction}

The $\CAT(0)$ comparison is a certain inequality for 6 distances between 4 points in a metric space.
The following descriptions, the so-called \emph{(2+2)-comparison}, is the most standard,
we refer to \cite{alexander-kapovitch-petrunin-2019,alexander-kapovitch-petrunin-2021} for other definitions and their equivalences.

Given a quadruple of points $p,q,x,y$ in a metric space $X$,
consider two \emph{model triangles}\footnote{that is, a plane triangle with the same sides}
$[\tilde p\tilde x\tilde y]=\tilde\triangle(pxy)$ 
and 
$[\tilde q\tilde x\tilde y]=\tilde\triangle(qxy)$ with common side $[\tilde x\tilde y]$.

\begin{wrapfigure}{r}{25mm}
\vskip-4mm
\centering
\includegraphics{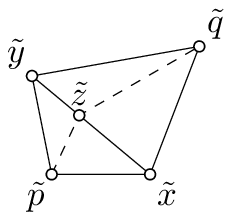}
\end{wrapfigure}

If the inequality
\[|p-q|_X\le |\tilde p-\tilde z|+|\tilde z-\tilde q|\]
holds for any point $\tilde z\in [\tilde x\tilde y]$, then we say that 
the quadruple $p,q,x,y$ \emph{satisfies $\CAT(0)$ comparison};
here $|p-q|_X$ denotes the distance from $p$ to $q$ in $X$.

If $\CAT(0)$ comparison holds for any quadruple (and any of its relabeling) in a metric space $X$,
then we say that \emph{$X$ is $\CAT(0)$}.

It is not hard to check that if a quadruple of points satisfies $\CAT(0)$ comparison for all relabeling,
then it admits a distance-preserving inclusion into a length $\CAT(0)$ space.
The following theorem generalizes this statement to 5-point metric spaces.

\begin{thm}{Toyoda's theorem}
Let $P$ be a 5-point metric space that satisfies $\CAT(0)$ comparison.
Then $P$ admits a distance-preserving inclusion into a length $\CAT(0)$ space $X$.

Moreover,
$X$ can be chosen to be a subcomplex of a 4-simplex such that (1) each simplex in $X$ has Euclidean metric and (2) the inclusion maps the 5 points on $P$ to the vertexes of the simplex.
\end{thm}

A slightly weaker version of this theorem was proved by Tetsu Toyoda \cite{toyoda}.
Our proof is shorter; it uses the fact that convex spacelike hypersurfaces in $\RR^{3,1}$ equipped with the induced length metrics are $\CAT(0)$ spaces \cite{milka}.
We construct a distance-preserving inclusion $\iota$ of $P$ into $\RR^4$ or $\RR^{3,1}$.
In the case of $\RR^4$ the convex hull $K$ of $\iota(P)$ can be taken as $X$;
in the case of $\RR^{3,1}$ we take as $X$ a spacelike part of the boundary of $K$.

It is expected that \emph{any 5 point metric space $P$ as in the theorem admits a distance-preserving inclusion in a product of trees}.\footnote{Later we found a counterexample \cite{petrunin-2021}.}

An analog of Toyoda's theorem does not hold for 6-point sets.
It can be seen by using the so-called $(4{+}2)$-comparison introduced in \cite{alexander-kapovitch-petrunin-2011};
this comparison holds for any length $\CAT(0)$ space, but may not hold for a space with $\CAT(0)$ comparison (if it is not a length space).

{\sloppy

The $(4{+}2)$-comparison is not a sufficient condition for $6$-point spaces.
More precisely, there are 6-point metric spaces that satisfy $(4{+}2)$ and $(2{+}2)$-comparisons but do not admit a distance-preserving embedding into a length $\CAT(0)$ space.
An example was constructed by the first author; it is described in \cite{alexander-kapovitch-petrunin-2011} right after 7.2.
See the final section for related questions.

}

\parbf{Acknowledgment.}
We want to thank Stephanie Alexander,
Yuri Burago,
and the anonymous referee for help.

The first author was partially supported by RFBR grant   20-01-00070, 
the second author was partially supported by NSF grant DMS-2005279.

\section{5-point arrays in 3-space}

Denote by $\mathcal{A}$ the space of all $5$ point arrays in $\RR^3$ that is nondegenerate in the following sense: (1) all 5 points do not lie on one plane and (2) no three points lie on one line.
Note that $\mathcal{A}$ is connected.

A $5$ point array  $x_1,\dots,x_5\in \RR^3$ defines an affine map from a 4-simplex to~$\RR^3$.
Fix an orientation of the 4-simplex and consider the induced orientations on its 5 facets.
Each facet may be mapped in an orientation-preserving, degenerate, or orientation-reversing way.
For each array consider the triple of integers $(n_+, n_0,n_-)$,
where $n_+$, $n_0$, and $n_-$ denote the number of orientation-preserving, degenerate, or orientation-reversing facets respectively.

Clearly $n_++n_0+n_-=5$ and since all 5 points cannot lie in one plane, we have that $n_+\ge 1$, $n_-\ge 1$, and $n_0\le1$.
Therefore, the value $m=n_- - n_+$ can take an integer value between $-3$ and $3$;
in this case, we say that an array belongs to $\mathcal{A}_m$.

It defines a subdivision of  $\mathcal{A}$ into 7 subsets $\mathcal{A}_{-3},\dots, \mathcal{A}_{3}$ with combinatorial configuration as on the diagram;
quadruples in one plane are marked in gray and the triple $(n_+, n_0,n_-)$ is written below.

\begin{figure}[ht!]
\centering
\begin{lpic}[t(-0mm),b(2mm),r(0mm),l(0mm)]{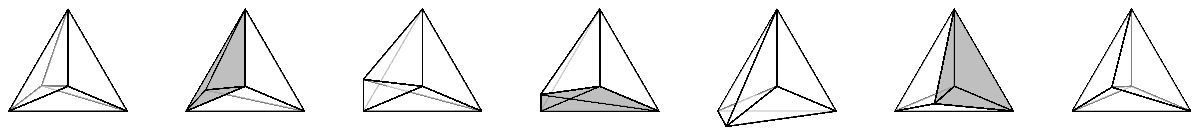}
\lbl[t]{6,0;$\mathcal{A}_{-3}$} 
\lbl[t]{6,-5;$(1,0,4)$} 
\lbl[t]{24,0;$\mathcal{A}_{-2}$}
\lbl[t]{24,-5;$(1,1,3)$} 
\lbl[t]{42,0;$\mathcal{A}_{-1}$} 
\lbl[t]{42,-5;$(2,0,3)$} 
\lbl[t]{60,0;$\mathcal{A}_{0}$} 
\lbl[t]{60,-5;$(2,1,2)$} 
\lbl[t]{78,0;$\mathcal{A}_{1}$}
\lbl[t]{78,-5;$(3,0,2)$} 
\lbl[t]{96,0;$\mathcal{A}_{2}$} 
\lbl[t]{96,-5;$(3,1,1)$} 
\lbl[t]{114,0;$\mathcal{A}_{3}$} 
\lbl[t]{114,-5;$(4,0,1)$} 
\end{lpic}
\vskip 4mm
\end{figure}

Every two quadrilaterals in the array have 3 common points that define a plane.
If the remaining two points lie on opposite sides from the plane,
then the corresponding facets have the same orientation;
if they lie on one side, then the orientations are opposite.
Therefore, the 7 subsets $\mathcal{A}_{-3},\dots \mathcal{A}_{3}$ can be described in the following way:

$\mathcal{A}_{-3}$ --- a tetrahedron with preserved orientation and one point inside.

$\mathcal{A}_{-2}$ --- a tetrahedron with preserved orientation and one point on a facet.

$\mathcal{A}_{-1}$ --- a double triangular pyramid formed by two tetrahedrons with preserved orientation.

$\mathcal{A}_{0}$ --- a pyramid over a convex quadrilateral 

$\mathcal{A}_{1}$ --- a double triangular pyramid formed by two tetrahedrons with reversed orientation.

$\mathcal{A}_{2}$ --- a tetrahedron with reversed orientation and one point on a facet.

$\mathcal{A}_{3}$ --- a tetrahedron with reversed orientation and one point inside.

Note that the complement $\spc{A}\backslash \spc{A}_0$ has two connected components formed by $\mathcal{A}_{-}=\mathcal{A}_{-3}\cup \mathcal{A}_{-2}\cup\mathcal{A}_{-1}$ and $\mathcal{A}_{+}=\mathcal{A}_{3}\cup \mathcal{A}_{2}\cup\mathcal{A}_{1}$.
Observe that each array in $\mathcal{A}_{-}$ has at least 3 positively oriented facets and each array in $\mathcal{A}_{+}$ has at least 3 negatively oriented facets.

\begin{thm}{Observation}\label{obs:connectedA}
Let $Q$ be a connected subset of $\spc{A}$ that does not intersect~$\spc{A}_0$.
Then either $Q\subset \spc{A}_+$ or $Q\subset \spc{A}_-$.
\end{thm}

\section{Associated form}

In this section we recall some facts about the so-called \emph{associated} form introduced in \cite{petrunin-2017};
it is a quadratic form 
$W_{\bm{x}}$ on $\RR^{n-1}$ associated
to a given $n$-point array $\bm{x}\z=(x_1,\dots,x_n)$ in a metric space $X$.

\parbf{Construction.}
Let $\triangle$ be the standard simplex $\triangle$ in $\RR^{n-1}$; that is, the first $(n-1)$ of its vertices $v_1,\dots,v_n$ form the standard basis on $\RR^{n-1}$,
 and $v_n=0$.

Recall that $|a-b|_X$ denotes the distance between points $a$ and $b$ in the metric space $X$.
Set
\[W_{\bm{x}}(v_i-v_j)=|x_i-x_j|^2_X\] 
for all $i$ and $j$.
Note that this identity defines $W_{\bm{x}}$ uniquely.

The constructed quadratic form $W_{\bm{x}}$ will be called the \emph{form associated to the point array $\bm{x}$}.

Note that an array $\bm{x}=(x_1,\dots,x_n)$ in a metric space $X$ is isometric to an array in Euclidean space if and only if 
$W_{\bm{x}}(v)\ge 0$
for any $v\in \RR^{n-1}$.

In particular, the condition $W_{\bm{x}}\ge 0$ for a triple $\bm{x}=(x_1,x_2,x_3)$ means that 
all three triangle inequalities for the distances between $x_1$, $x_2$, and $x_3$ hold.
For an $n$-point array, it implies that $W_{\bm{x}}(v)\ge 0$ for any vector $v$ in a plane spanned by a triple $v_i,v_j,v_k$.
In particular, we get the following:

\begin{thm}{Observation}\label{triangle-inq}
Let $W_{\bm{x}}$ be a form on $\RR^{n-1}$ associated with a point array $\bm{x}\z=(x_1,\dots,x_n)$.
Suppose that $L$ is a subspace of $\RR^{n-1}$ such that
$W_{\bm{x}}(v)< 0$ for any nonzero vector $v\in L$.
Then the projections of any 3 vertices of $\triangle$ to the quotient space $\RR^{n-1}/L$ are not collinear.
\end{thm}

\parbf{CAT(0) condition.}
Consider a point array $\bm{x}$ with 4 points.
From \ref{triangle-inq}, 
it follows that $W_{\bm{x}}$ 
is nonnegative on every plane parallel to a face of the tetrahedron $\triangle$.
In particular, $W_{\bm{x}}$ can have at most one negative eigenvalue.

Assume $W_{\bm{x}}(w)<0$ for some $w\in\RR^3$.
From \ref{triangle-inq}, the line $L_w$ spanned by
$w$ is transversal to each of 4 planes parallel to a face of $\triangle$.

Consider the projection of $\triangle$ along $L_w$ to a transversal plane. 
The projection of the 4 vertices of $\triangle$ lie in general position; 
that is, no three of them lie on one line.
Therefore, we can see one of two combinatorial pictures shown on the diagram.
Since the set of lines $L_w$ with $W_{\bm{x}}(w)<0$ is connected,
the combinatorics of the picture does not depend on the choice of $w$.

{

\begin{wrapfigure}{r}{33mm}
\vskip-0mm
\centering
\includegraphics{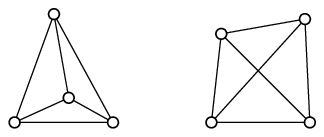}
\end{wrapfigure}

\begin{thm}{Claim}
If $\CAT(0)$ comparison holds in $X$, then the diagram on the right cannot appear. 
\end{thm}

(The converse holds as well, but we will not need it.)

}

\parit{Proof.}
Suppose we see the picture on the right.

Let $[v_1,v_3]$ and $[v_2,v_4]$ be the line segments of $\triangle$ that correspond to the diagonals on the picture.
Denote by $m$ the point of $[v_1,v_3]$ that corresponds to the point of intersection.

\begin{wrapfigure}{r}{33mm}
\vskip-4mm
\centering
\includegraphics{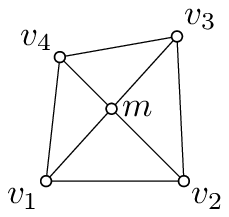}
\end{wrapfigure}

In the plane spanned by $[v_2,v_4]$ and $w$, the vector $w$ is timelike.
Therefore we have the following reversed triangle inequality:
\[|v_2-m|+|v_4-m|<|v_2-v_4|;\]
here we use shortcut $|a-b|=\sqrt{W(a-b)}$.

Note that the triangles $[v_1v_2v_3]$ and $[v_1v_3v_4]$ with metric induced by $W$ are isometric to model triangles of $[x_1x_2x_3]$ and $[x_1x_3x_4]$.
Whence (2+2)-point comparison does not hold.
\qeds

The claim implies the following:

\begin{thm}{Observation}\label{cat0-proj}
Suppose a metric on $\bm{x}\z=(x_1,\dots,x_n)$ satisfies $\CAT(0)$ comparison
and $W_{\bm{x}}$ is its associated form on $\RR^{n-1}$.
Assume that $L$ is a subspace of $\RR^{n-1}$ such that
$W_{\bm{x}}(v)< 0$ for any nonzero vector $v\in L$.
Then if the projections of 4 vertices of $\triangle$ to the quotient space $\RR^{n-1}/L$ lies in one plane, then its projection looks like the picture on the left;
that is, one of the points lies in the triangle formed by the remaining three points.
\end{thm}

\begin{thm}{Corollary}\label{cor:3+2}
Suppose a metric on $\bm{x}=(x_1,\dots,x_5)$ satisfies $\CAT(0)$ comparison
and $W_{\bm{x}}$ is its associated form on $\RR^{4}$.
Assume that $L$ is a subspace of $\RR^{4}$ such that
$W_{\bm{x}}(v)< 0$ for any nonzero vector $v\in L$.
Then $\dim L\le 1$.

Moreover, if $\dim L= 1$, then the projections of the vertices of $\triangle$ to the quotient space $\RR^3=\RR^4/L$ belong to $\spc{A}\backslash\spc{A}_0$ (defined in the previous section). 
\end{thm}

\parit{Proof.}
If $\dim L\ge 2$, then $\dim (\RR^4/L)\le 2$.
By \ref{triangle-inq}, these 5 projections lie in a general position; that is, no three of these projections lie on one line. 
Therefore, $\RR^4/L=2$ is the plane.

Any 5 points in a general position on the plane include 4 vertices of a convex quadrangle.
The latter contradicts \ref{cat0-proj}.
\qeds

\section{Convex spacelike surfaces}

\begin{wrapfigure}{r}{55mm}
\vskip-4mm
\centering
\includegraphics{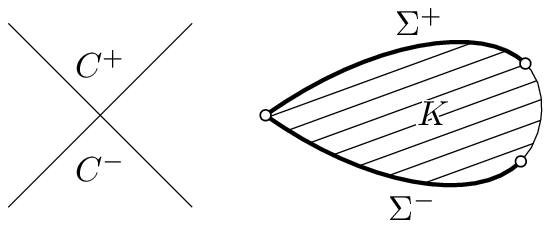}
\end{wrapfigure}

Let $W$ be a quadratic form on $\RR^4$.
Suppose that  $W$ has exactly one negative eigenvalue.
Choose future and past cones $C^+$ and $C^-$ for $W$;
that is, $C^+$ and $C^-$ are connected components of the set
$\set{v\in\RR^4}{W(v)< 0}$.
A subset $S$ in $\RR^4$ will be called \emph{spacelike} if $W(x-y)\ge 0$ for any $x,y\in S$.

Let $K$ be a convex body in $\RR^4$;
denote by $\Sigma$ the surface of $K$. 
A point $p$ lies on the \emph{upper side} of $\Sigma$ (briefly $p\in\Sigma^+$) if there is a spacelike hyperplane in $\RR^4$ that supports $\Sigma$ at $p$ from above;
more precisely if the Minkowski sum $\{p\}+C^+$ does not intersect~$K$.

Similarly, we define the \emph{lower side} of $\Sigma$ denoted by $\Sigma^-$.
Note that $\Sigma^+$ and $\Sigma^-$ might have common points.
The subsets $\Sigma^+$ and $\Sigma^-$ are spacelike;
in particular, the length of any Lipschitz curve in these subsets can be defined and it leads to induced intrinsic pseudometrics on $\Sigma^+$ and $\Sigma^-$.
Abusing notation, we will not distinguish a pseudometric space and the corresponding metric space. 

\begin{thm}{Lemma}\label{lem:sides}
Let $\Sigma$ be the surface of a convex set $K$ in $\RR^4$ and $C^\pm$ be the future and past cones for a quadratic form $W$.
Then the upper and lower sides $\Sigma^+$ and $\Sigma^-$ of $\Sigma$ equipped with the induced intrinsic metric are $\CAT(0)$ length spaces.

Moreover, if a \emph{line segment} $[pq]$ in $\RR^4$ 
lies on $\Sigma^\pm$, then $[pq]$ is a minimizing geodesic in $\Sigma^\pm$; that is,
\[|p-q|_{\Sigma^\pm}^2=W(p-q).\]

\end{thm}

This lemma is essentially stated by Anatolii Milka \cite[Theorem~4]{milka}; we give a sketch of alternative proof based on smooth approximation.

\parit{Sketch.}
We can assume that $W$ is nondegenerate; that is, after a linear change of coordinates it is the standard form on $\RR^{3,1}$.
If not, then there is a $W$-preserving projection of $\RR^4$ to a $W$-nondegenerate subspace; apply this projection and note that this subspace is isometric a subspace of $\RR^{3,1}$.

Assume $S$ is a smooth strictly spacelike hypersurface in $\RR^{3,1}$ with convex epigraph.
By Gauss formula, $S$ has nonpositive sectional curvature.

Suppose a strictly spacelike hyperplane $\Pi$ cuts from $S$ a disc $D$.
Recall that Liberman's lemma \cite[Theorem~3]{milka} implies that time coordinate is convex on any geodesic in $S$.
We may assume that time is vanishing on~$\Pi$;
therefore, by the lemma, $D$ has a convex set in $S$.
Therefore the Cartan--Hadamard theorem \cite{alexander-kapovitch-petrunin-2021} implies that that $D$ is $\CAT(0)$.

Now suppose $D_n$ is a sequence of smooth discs of the described type that converges to a (possibly nonsmooth) disc $D$.
Note that the metric on $D_n$ converges to the induced pseudometric on $D$.
It follows that the metric space $D'$ that corresponds to $D$ is $\CAT(0)$.

The disc $D$ might contain lightlike segments which have zero length.
Note that every maximal lightlike segment in $D$ starts at its interior point and goes to the boundary.
Consider the map $\iota\:D\to D$ that sends each maximal lightlike segment to its starting point.
Note that the sublemma below implies that $\iota$ is length-nonincreasing.
Since $|x-\iota(x)|_D=0$, we get that the $D'$ is isometric to the image of $\iota$ with the induced metric.

\begin{wrapfigure}{r}{61mm}
\vskip-0mm
\centering
\includegraphics{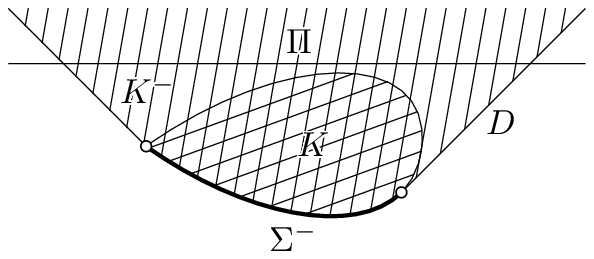}
\end{wrapfigure}

Consider the Minkowski sum
\[K^-=K+C^+;\]
it has a convex spacelike boundary $\partial K^-$.
Choose a strictly spacelike hyperplane $\Pi$ that lies above $K$.
Denote by $D$ the subset of $\partial K^-$ below $\Pi$.
Let us equip $D$ with induced intrinsic pseudometric.
By construction $\Sigma^-$ is isometric to $\iota(D)$.
It follows that $\Sigma^-$ is $\CAT(0)$. 

Now suppose a line segment $[pq]$ in $\RR^4$ 
lies on $\Sigma^-$.
Choose a supporting hyperplane $\Pi$ at the midpoint of $[pq]$.
Choose time coordinate that vanish on $\Pi$;
by Liberman's lemma, every shortest path in $\Sigma^-$ between $p$ and $q$ has to lie on $\Pi$;
that is, the intersection $\Sigma^-\cap \Pi$ is a convex subset of $\Sigma^-$.
Therefore $[pq]$ is convex in $\Sigma^-$ which implies the second statement.
\qeds

\begin{thm}{Sublemma}
Let $u$ and $v$ be two lightlike vectors in $\RR^{3,1}$.
Suppose that the union of two half-lines $s\mapsto p+s\cdot u$ and $t\mapsto q+t\cdot v$ for $s,t\ge 0$ is a spacelike set.
Then the function 
$(s,t)\mapsto |(p+s\cdot u)-(q+t\cdot v)|$
is nondecreasing in both arguments,
where $|w|\df\sqrt{\langle w,w\rangle}$ for a spacelike vector~$w$.
\end{thm}

\parit{Proof.}
Since $u$ and $v$ are lightlike, $\langle u,u\rangle=\langle v,v\rangle=0$.
Since the union of two half-lines is spacelike, $(p+s\cdot u)-(q+t\cdot v)$ is spacelike for any $s,t\ge 0$.
It follows that
\begin{align*}
0&\le |(p+s\cdot u)-(q+t\cdot v)|^2=
\\
&=|p-q|^2-2\cdot s \cdot\langle u,q-p\rangle-2\cdot t\cdot\langle v,p-q\rangle- 2\cdot s\cdot t\cdot \langle u,v\rangle
\end{align*}
for any $s,t\ge0$.
Therefore
\begin{align*}
\langle u,q-p\rangle&\le 0,
& 
\langle v,p-q\rangle&\le 0
&
\langle u,v\rangle&\le 0.
\end{align*}
Whence the result.
\qeds

Assume $v$ is a nonzero vector in $\RR^4$ and $p\in\Sigma$.
We say that $p$ lies on the \emph{upper side of $\Sigma$ with respect to $v$} (briefly $p\in \Sigma^+(v)$) if $p+t\cdot v\notin K$ for any $t>0$.
Correspondingly, $p$ lies on the \emph{lower side of $\Sigma$ with respect to $v$} (briefly $p\in \Sigma^-(v)$) if $p+t\cdot v\notin K$ for any $t<0$.

\begin{thm}{Observation}\label{obs:Sigma(v)}
Let $K$ be a compact convex set in $\RR^4$ and $C^\pm$ be the future and past cones for a quadratic form $W$.
Then the upper (lower) side of the boundary surface $\Sigma$ of $K$ can be described as the intersection of the upper (respectively lower) sides of $\Sigma$ with respect to all vectors $v\in C^+$;
that is,
\[\Sigma^\pm=\bigcap_{v\in C^+}\Sigma^\pm(v).\]
\end{thm}

\section{Proof assembling}

\parit{Proof of Toyoda's theorem.}
Let $\{x_1,\dots,x_5\}$ be the points in $P$.
Choose a 5-simplex $\triangle$ in $\RR^4$; denote by $W$ the form associated with the point array $(x_1,\z\dots,x_5)$.

If $W\ge0$, then $P$ admits a distance preserving embedding into Euclidean 4-space, so one can take the convex hull of its image as $X$.

Suppose $W(v)<0$ for some $v\in\RR^4$.
Since $P$ is $\CAT(0)$, \ref{cor:3+2} implies that $W$ has exactly one negative eigenvalue.
Moreover, if a line $L$ is spanned by a vector $v$ such that $W(v)<0$, then the projection of the vertices of the simplex to $\RR^3=\RR^4/L$ belongs to $\spc{A}\backslash\spc{A}_0$.

The space of such lines $L$ is connected.
By \ref{obs:connectedA}, we can assume that all the projections belong to $\spc{A}_-$.
That is, we can choose timelike orientation such that for any $v\in C^+$ the lower part $\Sigma^-(v)$ of $\Sigma=\partial \triangle$ has at least 3 facets of $\triangle$.

In particular, $\Sigma^-(v)$ contains all edges of $\triangle$ for any $v\in C^+$.
By \ref{obs:Sigma(v)}, $\Sigma^-$ contains all edges of $\triangle$.
By \ref{lem:sides}, $\Sigma^-$ with induced (pseudo)metric is a length $\CAT(0)$ space.

Since all edges of $\triangle$ lie in $\Sigma^-$, the inclusion $P\hookrightarrow \Sigma^-$ is distance preserving.
Whence we can take $X=\Sigma^-$.

Finally, observe that in each case $X$ is a subcomplex of $\triangle$ that includes all edges and has a model metric on each simplex.
\qeds

\section{Remarks}

Let us recall the definition of \emph{graph comparison} given by Vladimir Zolotov and the authors \cite{lebedeva-petrunin-zolotov} and use it to formulate a few related questions.

Let $\Gamma$ be a graph with vertices $v_1,\dots,v_n$.
A metric space $X$ is said to meet the $\Gamma$-comparison if for any set of points in $X$ labeled by vertices of $\Gamma$ there is a model configuration $\tilde v_1,\dots,\tilde v_n$ in the Hilbert space $\HH$ such that 
if $v_j$ is adjacent to $v_j$, then
\[|\tilde v_i-\tilde v_j|_{\HH}\le | v_i-v_j|_{X}\]
and
if $v_j$ is nonadjacent to $v_j$, then
\[|\tilde v_i-\tilde v_j|_{\HH}\ge | v_i-v_j|_{X}.\]

The $C_4$-comparison (for the 4-cycle $C_4$ on the diagram) defines $\CAT(0)$ comparison.
Tetsu Toyoda have shown that $C_4$-comparison imlies graph comparisons for all cycles $C_n$ \cite{toyoda2020}; remakably, the metric space is \emph{not} assumed to be intrinsic. 
\begin{figure}[ht!]
\vskip-0mm
\centering
\includegraphics{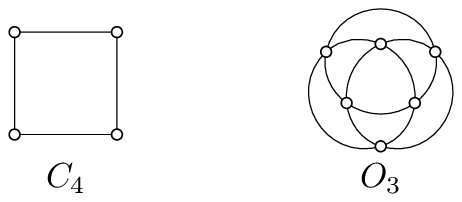}
\end{figure}
The $O_3$-comparison (for the octahedron graph $O_3$ on the diagram) defines another comparison.
Since $O_3$ contains $C_4$ as an induced subgraph, we get that $O_3$-comparison is stronger than $C_4$-comparison.

\begin{thm}{Open question}
Is it true that octahedron-comparison holds in any 6 points in a length $\CAT(0)$ space?

And, assuming the answer is affirmative, what about the converse: is it true that any 6-point metric space that satisfies octahedron-comparison admits a distance preserving embedding in a length $\CAT(0)$ space?
\end{thm}

The analogous questions for  spaces with nonnegative curvature in the sense of Alexandrov (briefly $\CBB(0)$) are open as well.
The $\CBB(0)$ comparison is equivalent to the $3$-tree comparison (for the tripod-tree shown first on the following diagram).
\begin{figure}[h!]
\vskip-0mm
\centering
\includegraphics{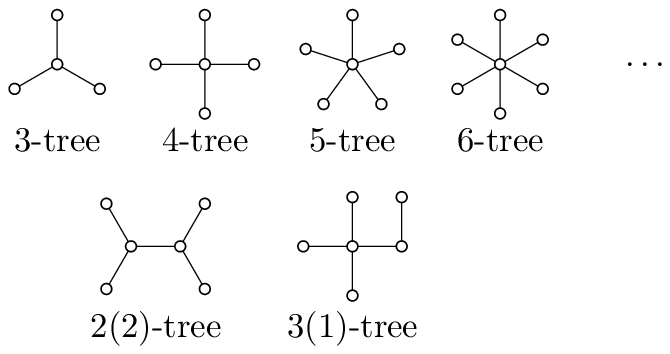}
\end{figure}
It turns out that any length $\CBB(0)$ space satisfies the comparison for the other trees on the diagram; it is formed by an infinite family of star-shaped trees and two trees with 6 vertices \cite{alexander-kapovitch-petrunin-2011,lebedeva-petrunin-zolotov}.
(The 4-tree comparison (the second tree on the diagram) is equivalent to the so-called $(4{+}1)$-point comparison in the terminology of \cite{alexander-kapovitch-petrunin-2011}.)

We expect that this comparison provides a necessary and sufficient condition for 5-point sets.
Namely, we expect an affirmative answer to the following stronger question.

\begin{thm}{Question}
Suppose a 5-point metric space $P$ satisfies the $4$-tree comparison.
Is it true that $P$ admits a distance preserving embedding into a length $\CBB(0)$ space?
\end{thm}

Finally, let us mention a related question about a 6-point condition. 

\begin{thm}{Question}
Suppose a 6-point metric space $P$ satisfies the 5-tree, 2(2)-tree, and 3(1)-tree comparisons.
Is it true that $P$ admits a distance preserving embedding into a length $\CBB(0)$ space?
\end{thm}


{\sloppy
\printbibliography[heading=bibintoc]
\fussy
}

\Addresses
\end{document}